\documentclass{article}
\usepackage{amssymb}
\usepackage{amsmath}
\usepackage{abstract}
\usepackage[utf8]{inputenc}
\usepackage[hidelinks]{hyperref}

\usepackage{titlesec}
\titleformat{\section}[block]{\center\scshape\large}{\thesection.}{0.7em}{}
\titleformat{\subsection}[block]{\center}{\thesubsection.}{0.7em}{}
\renewcommand\thesection{\Roman{section}}

\usepackage{amsthm}
\newtheorem{theorem}{Theorem}
\newtheorem{lemma}{Lemma}

\theoremstyle{definition}

\newcommand{\keywords}[1]
{
	\small
	\textbf{Keywords:} #1
}

\title{A Simple Proof of Siegel's Theorem Using Mellin Transform}

\usepackage{authblk}
\author[a]{Zihao Liu}
\affil[a]{International Department, The Affiliated High School of SCNU,\authorcr Email: \url{mailto:travor_lzh@163.com}}

\date{}

\begin{document}
	\maketitle
	\begin{abstract}
		In this paper, we present a simple analytic proof of Siegel's theorem that concerns the lower bound of $L(1,\chi)$ for primitive quadratic $\chi$. Our new method compares an elementary lower bound with an analytic upper bound obtained by the inverse Mellin transform of $\Gamma(s)$.
	\end{abstract}
	\section{Introduction}
	In 1935, Siegel\cite{Siegel1935} introduces the function
	\begin{equation}
		\label{eqnfs}
		f(s)=\zeta(s)L(s,\chi_1)L(s,\chi_2)L(s,\chi_1\chi_2)
	\end{equation}
	where $\chi_1$ and $\chi_2$ are primitive quadratic characters modulo $q_1$ and $q_2$ respectively. By exploring its algebraic properties, he shows that a very strong lower bound can be established for $L(1,\chi)$:
	\begin{theorem}[Siegel]
		\label{thsiegel}
		For all $\varepsilon>0$ there exists a constant $C(\varepsilon)>0$ such that
		\begin{equation}
			L(1,\chi)>C(\varepsilon)q^{-\varepsilon}
		\end{equation}
		holds for any primitive quadratic $\chi$ modulo $q$.
	\end{theorem}
	Although the statement of \autoref{thsiegel} is analytic, it leads to strong conclusions in the distribution of prime numbers in arithmetic progressions. Using this result, Walfisz\cite{walfisz_zur_1936} improved the zero-free region of $L(s,\chi)$ to obtain the prime number theorem for arithmetic progressions in the following form:
	\begin{theorem}[Siegel-Walfisz]
		\label{thsw}
		Let $\pi(x;q,a)$ denotes the number of primes $\le x$ that are $\equiv a\pmod q$. Then for all $A>0$, there exists $C_A>0$ such that the following estimate
		\begin{equation}
			\label{eqsw}
			\pi(x;q,a)={1\over\varphi(q)}\int_2^x{\mathrm dt\over\log t}+O_A\left\{x\exp\left(-C_A\sqrt{\log x}\right)\right\}
		\end{equation}
		holds when $(a,q)=1$ and $q\le(\log x)^A$.
	\end{theorem}
	This result is very powerful in additive problems concerning primes. For instance, Vinogradov deduces from \autoref{thsw} that every large odd integer is a sum of three primes.\footnote{This is often known as the ternary Goldbach's conjecture. See §25 and §26 of \cite{davenport_multiplicative_1980} for an account of Vinogradov's proof.} Under \autoref{thsw}, Mirsky\cite{mirsky_number_1949} shows that every large integer is a sum of a prime and a $k$-free integer.
	
	The original proof of \autoref{thsiegel} uses algebraic number theory. Later in 1949, Estermann\cite{Estermann1949} obtained a simple proof using purely analytic methods. Few decades after that, Goldfeld\cite{goldfeld_simple_1974}\cite{liu2021goldfelds} gave a much more simplified analytic the proof using contour integration. In this paper, we propose a new contour-integration proof of \autoref{thsiegel} based on the inverse Mellin transform of $\Gamma(s)$. In particular, the new approach simplifies Goldfeld's method because it uses Abelian summation instead of Ces\`aro summation\footnote{A detailed description of these summation methods is accessible in §5.2 of \cite{montgomery_multiplicative_2007}}.
	\section{Analytic Lemmas}
	From now on, $s=\sigma+it$ always denotes a complex number with an abscissa of $\sigma$ and an ordinate of $t$.
	\begin{lemma}[Phragm\'en-Lindel\"of]
		\label{lmpl}
		If $\phi(s)$ is regular and $O_\varepsilon(e^{\varepsilon|t|})$, for any $\varepsilon>0$, in the strip $\sigma_1\le\sigma\le\sigma_2$, and
		$$
		\phi(\sigma_1+it)\ll|t|^{k_1},\quad\phi(\sigma_2+it)\ll|t|^{k_2}
		$$
		then $\phi(s)\ll|t|^{\max(k_1,k_2)}$ holds uniformly in $\sigma_1\le\sigma\le\sigma_2$.
	\end{lemma}
	\begin{proof}
		See §5.65 of \cite{titchmarsh_theory_2002}.
	\end{proof}
	\begin{lemma}
		\label{lmfb}
		Let $f(s)$ be defined as in \eqref{eqnfs}, then for all $\varepsilon>0$. The estimate
		\begin{equation}
			\label{eqnfb}
			f(s)\ll_\varepsilon(q_1q_2)^{1+\varepsilon}|t|^{2+\varepsilon}
		\end{equation}
		holds uniformly in $\sigma\ge0$.
	\end{lemma}
	\begin{proof}
		It is well known that when $\sigma$ lies in a fixed strip and $|t|\to\infty$, $\zeta(s)$ and $L(s,\chi)$ satisfies the following asymptotic functional equations\footnote{See §10.1 of \cite{montgomery_multiplicative_2007} for a full derivation.}:
		\begin{equation}
			\label{eqnzf}
			\zeta(s)\ll|t|^{1/2-\sigma}|\zeta(1-s)|
		\end{equation}
		\begin{equation}
			\label{eqnlf}
			L(s,\chi)\ll(q|t|)^{1/2-\sigma}|L(1-s,\overline\chi)|
		\end{equation}
		where $\chi$ is a primitive character modulo $q$.
		for . Since $\zeta(s)$ and $L(s,\chi)$ converge absolutely for all $\sigma\ge1+\varepsilon$, we see that when $\sigma=-\varepsilon$, \eqref{eqnzf} and \eqref{eqnlf} can be simplified into
		\begin{equation}
			\label{eqnzf2}
			\zeta(s)\ll_\varepsilon|t|^{1/2+\varepsilon}
		\end{equation}
		\begin{equation}
			\label{eqnlf2}
			L(s,\chi)\ll_\varepsilon(q|t|)^{1/2+\varepsilon}
		\end{equation}
		Plugging \eqref{eqnzf2} and \eqref{eqnlf2} into \eqref{eqnfs}, we see that \eqref{eqnfb} holds for $\sigma=-\varepsilon/4$, and finally we can apply \autoref{lmpl} to extend this estimate to $\sigma\ge-\varepsilon/4$.
	\end{proof}
	\begin{lemma}
		Let $\Gamma(s)$ denote the Gamma function. Then the estimate
		\begin{equation}
			\Gamma(s)\ll|t|^{\sigma-1/2}e^{-\pi|t|/2}
		\end{equation}
		holds whenever $\sigma$ lies in a fixed interval.
	\end{lemma}
	\begin{proof}
		By Stirling's formula, we know that when $\sigma$ lies in a fixed interval
		\begin{align*}
			\log\Gamma(s)
			&=\left(s-\frac12\right)\log s-s+O(1) \\
			&=\left(s-\frac12\right)\log(it)+\left(s-\frac12\right)\log\left(1+{\sigma\over it}\right)-s+O(1) \\
			&=\left(s-\frac12\right)\log(it)+it\log\left(1+{\sigma\over it}\right)-s+O(1) \\
			&=\left(s-\frac12\right)\log(it)-it+O(1)
		\end{align*}
		Taking real parts on both side, we see that as $t\to+\infty$ there is
		\begin{equation}
			\label{eqnstirling}
			\log|\Gamma(s)|=\left(\sigma-\frac12\right)\log t-{\pi t\over2}+O(1)
		\end{equation}
		Therefore, exponentiating on both side of \eqref{eqnstirling} yields the desired result.
	\end{proof}
	\begin{lemma}\label{lmgamma}
		For all $y>0$ there is
		\begin{equation}
			e^{-y}={1\over2\pi i}\int_{2-i\infty}^{2+i\infty}y^{-s}\Gamma(s)\mathrm ds
		\end{equation}
	\end{lemma}
	\begin{proof}
		The result follows directly by applying Mellin's inversion formula to
		\begin{equation}
			\Gamma(s)=\int_0^\infty y^{s-1}e^{-y}\mathrm dy
		\end{equation}
	\end{proof}
	\begin{lemma}
		\label{lmz}
		For $0<\sigma<1$, we have $\zeta(\sigma)<0$.
	\end{lemma}
	\begin{proof}
		By partial summation, we have
		\begin{equation}
			\zeta(\sigma)={\sigma\over\sigma-1}-\sigma\int_1^\infty{\{x\}\over x^{\sigma+1}}\mathrm dx\le{\sigma\over\sigma-1}
		\end{equation}
		The right hand side immediately concludes the proof.
	\end{proof}
	\begin{lemma}
		\label{lml1}
		Let $\chi$ be nonprincipal character modulo $q$ then $L(1,\chi)\ll\log q$.
	\end{lemma}
	\begin{proof}
		Using the fact that $|\chi(n)|\le1$, we have
		\begin{align*}
			|L(1,\chi)|
			&\le\left|\sum_{n\le N}{\chi(n)\over n}\right|+\left|\sum_{n>N}{\chi(n)\over n}\right| \\
			&\le\sum_{n\le N}\frac1n+\left|\left[\frac1x\sum_{N<n\le x}\chi(n)\right]_N^\infty-\int_N^\infty\left(\sum_{N<n\le t}\chi(n)\right){\mathrm dt\over t^2}\right| \\
			&\le1+\sum_{2\le n\le N}\int_{n-1}^n{\mathrm dt\over t}+\max_x\left|\sum_{N<n\le x}\chi(n)\right|\left\{\frac1N+\int_N^\infty{\mathrm dt\over t^2}\right\} \\
			&\le1+\log N+{2\varphi(q)\over N}\ll\log N+\frac qN
		\end{align*}
		Setting $N=q$ gives the desired result.
	\end{proof}
	\begin{lemma}
		\label{lml}
		Let $\chi$ be a quadratic character such that $L(s,\chi)$ is free of real zeros in $s>1-\varepsilon$. Then $L(\beta,\chi)>0$ holds for any $1-\varepsilon<\beta<1$.
	\end{lemma}
	\begin{proof}
		Since $L(s,\chi)$ is continuous in $[1-\varepsilon,1]$ and $L(1,\chi)>0$\footnote{This is an auxiliary result used to prove Dirichlet's theorem. See §4.3 of \cite{montgomery_multiplicative_2007} for a derivation}, the result immediately follows.
	\end{proof}
	\begin{lemma}
		\label{lmbeta}
		For any $\varepsilon>0$ there exists a primitive quadratic $\chi_1$ modulo $q_1$ and $1-\varepsilon<\beta<1$ such that $f(\beta)\le0$ holds for all quadratic $\chi_2$ modulo $q_2$.
	\end{lemma}
	\begin{proof}
		On one hand, if no $\chi$ can be found such that $L(s,\chi)$ has a zero in $(1-\varepsilon,1)$. Then it follows from \autoref{lmz} and \autoref{lml} that $f(\beta)<0$ for all $1-\varepsilon<\beta<1$.

		On the other hand, if we are unable to find a quadratic primitive $\chi$ such that $L(s,\chi)$ does possess a real zero in $(1-\varepsilon,1)$. Then let $\chi_1=\chi$ and $\beta$ be the real zero so that $f(\beta)=0$. Consequently for every $\varepsilon>0$, there exists a primitive quadratic $\chi_1$ modulo $q_1$ and $1-\varepsilon<\beta<1$ such that $f(\beta)\le0$.
	\end{proof}
	\section{Proof of Siegel's Theorem}
	Similar to Goldfeld's method\cite{goldfeld_simple_1974}, our approach also studies the partial sum of $f(s)$
	\begin{equation}
		A(x,\beta)=\sum_{n\le x}{a_n\over n^\beta}
	\end{equation}
	where $a_n$ denote the Dirichlet series coefficient of $f(s)$ and $1-\varepsilon<\beta<1$ It follows from literature\footnote{See §21 of \cite{davenport_multiplicative_1980} for a detailed account} that $a_n\ge0$ and $a_1=1$, so we have $A(x,\beta)\ge1$ when $x\ge1$. In addition, because the exponential decay function satisfies
	\begin{equation*}
		e^{-n/x}
		\begin{cases}
			\ge0 & n>x \\
			\ge e^{-1} & n\le x
		\end{cases}
	\end{equation*}
	
	we also have
	\begin{equation}
		\label{eqnax}
		1\le A(x,\beta)\le e\sum_{n\ge1}{a_n\over n^\beta}e^{-x/n}
	\end{equation}
	Now, we apply \autoref{lmgamma} to the exponential function in \eqref{eqnax} so that
	\begin{equation}
		e^{-1}\le{1\over2\pi i}\int_{2-i\infty}^{2+i\infty}x^s\Gamma(s)f(s+\beta)\mathrm ds\triangleq I
	\end{equation}
	To estimate the integral, we move the path of integration to $\sigma=-\beta$ so that it follows from \autoref{lmfb} that
	\begin{align}
		I
		&=x^{1-\beta}\Gamma(1-\beta)\lambda+f(\beta)+x^{-\beta}\int_{-\beta-i\infty}^{-\beta+i\infty}|\Gamma(s)f(s+\beta)|\mathrm ds \\
		\label{eqni}&=x^{1-\beta}\Gamma(1-\beta)\lambda+f(\beta)+O_\varepsilon\left\{x^{-\beta}(q_1q_2)^{1+\varepsilon}\int_0^\infty t^{2+\varepsilon}|\Gamma(-\beta+it)|\mathrm dt\right\}
	\end{align}
	where $\lambda=L(1,\chi_1)L(1,\chi_2)L(1,\chi_1\chi_2)$ is the residue of $f(s)$ at $s=1$. Since $\Gamma(s+1)=s\Gamma(s)$, the remaining integral will be bounded by
	\begin{align*}
		\int_0^\infty t^{2+\varepsilon}|\Gamma(-\beta+it)|
		&\ll{1\over1-\beta}\int_0^\infty t^{5/2-\beta+\varepsilon}e^{-\pi t/2}\mathrm dt \\
		&\ll{1\over1-\beta}\int_0^\infty t^{3/2+\varepsilon}e^{-\pi t/2}\mathrm dt\ll_\varepsilon{1\over1-\beta}
	\end{align*}
	Now, if we choose $\beta$ and $\chi_1$ according to \autoref{lmbeta}, then we can ignore the $f(\beta)$ term to simplify \eqref{eqni} into
	\begin{equation}
		\label{eqni2}
		1\ll_\varepsilon x^{1-\beta}\lambda+x^{-\beta}q_2^{1+\varepsilon}
	\end{equation}
	in which all $\beta$ and $q_1$ terms in the coefficients are absorbed into $\ll_\varepsilon$. To simplify this even further, we set $x^\beta=cq_2^{1+\varepsilon}$ for some small $c>0$ so that the left hand side of \eqref{eqni2} will still be positive even after subtracted by $x^{-\beta}q_2^{1+\varepsilon}$. To further simplify the right hand side, we apply \autoref{lml1} to $\lambda$ so that for $q_2>q_1(\varepsilon)$ there is
	\begin{equation}
		1\ll_\varepsilon x^{1-\beta}\lambda\ll x^{1-\beta}(\log q_1)(\log q_1q_2)L(1,\chi_2)\ll_\varepsilon x^\varepsilon(\log q_2)L(1,\chi_2)
	\end{equation}
	Transforming this equation, we have
	\begin{align*}
		L(1,\chi_2)
		&\gg_\varepsilon x^{-\varepsilon}(\log q_2)^{-1} \\
		&\gg_\varepsilon q_2^{-\varepsilon(1+\varepsilon)/\beta}(\log q_2)^{-1} \\
		&\gg_\varepsilon q_2^{-\varepsilon{1+\varepsilon\over1-\varepsilon}}(\log q_2)^{-1}
	\end{align*}
	Without loss of generality, we assume $\varepsilon\le 1/2$, so that
	\begin{equation}
		L(1,\chi_2)\gg_\varepsilon q^{-3\varepsilon}(\log q_2)^{-1}\gg_\varepsilon q_2^{-4\varepsilon}
	\end{equation}
	This lower bound becomes \autoref{thsiegel} after a change of variable.
	\bibliographystyle{plain}
	\bibliography{refs.bib}

\begin{thebibliography}{1}

\bibitem{davenport_multiplicative_1980}
Harold Davenport.
\newblock {\em Multiplicative {Number} {Theory}}, volume~74 of {\em Graduate
  {Texts} in {Mathematics}}.
\newblock Springer New York, New York, NY, 1980.

\bibitem{Estermann1949}
T.~{Estermann}.
\newblock {On Dirichlet's \(L\) functions}.
\newblock {\em {J. Lond. Math. Soc.}}, 23:275--279, 1949.

\bibitem{goldfeld_simple_1974}
D.~M. Goldfeld.
\newblock A {Simple} {Proof} of {Siegel}'s {Theorem}.
\newblock {\em Proceedings of the National Academy of Sciences},
  71(4):1055--1055, April 1974.

\bibitem{liu2021goldfelds}
Zihao Liu.
\newblock On {Goldfeld's} {Proof} of {Siegel's} {Theorem}, 2021.
\newblock \url{https://arxiv.org/abs/2201.11145v1}.

\bibitem{mirsky_number_1949}
L.~Mirsky.
\newblock The {Number} of {Representations} of an {Integer} as the {Sum} of a
  {Prime} and a k-{Free} {Integer}.
\newblock {\em The American Mathematical Monthly}, 56(1):17, January 1949.

\bibitem{montgomery_multiplicative_2007}
Hugh~L. Montgomery and Robert~C. Vaughan.
\newblock {\em Multiplicative number theory {I}: classical theory}.
\newblock Number~97 in Cambridge studies in advanced mathematics. Cambridge
  University Press, Cambridge, UK ; New York, 2007.
\newblock OCLC: ocm61757122.

\bibitem{Siegel1935}
Carl Siegel.
\newblock Über die classenzahl quadratischer zahlkörper.
\newblock {\em Acta Arithmetica}, 1(1):83--86, 1935.

\bibitem{titchmarsh_theory_2002}
E.~C. Titchmarsh.
\newblock {\em The theory of functions}.
\newblock Oxford science publications. Oxford Univ. Press, Oxford, 2. ed.,
  reprinted edition, 2002.
\newblock OCLC: 249703508.

\bibitem{walfisz_zur_1936}
Arnold Walfisz.
\newblock Zur additiven {Zahlentheorie}. {II}.
\newblock {\em Math Z}, 40(1):592--607, December 1936.

\end{thebibliography}

\end{document}